\documentstyle[11pt,amssymb,amstex]{amsart}

\newtheorem{thm}{Theorem}[section]

\newtheorem{cor}[thm]{Corollary}
\newtheorem{prop}[thm]{Proposition}

\newtheorem{prob}[thm]{Problem}

\begin{document}

\title[Doubly stochastic quadratic operators and Birkhoff's problem]
{Doubly stochastic quadratic operators and Birkhoff's problem}
\author{Rasul Ganikhodzhaev}
\address{Rasul Ganikhodzhaev \\
Department of Mechanics and Mathematics,\\
National University of Uzbekistan\\
Vuzgorodok, Tashkent, 100174, Uzbekistan} \email{{\tt
rganikhodzhaev@@gmail.com}}
\author{Farruh Shahidi}
\address{Farruh Shahidi \\
Department of Mechanics and Mathematics,\\
National University of Uzbekistan\\
Vuzgorodok, Tashkent, 100174, Uzbekistan} \email{{\tt
farruh.shahidi@@gmail.com}}

\begin{abstract}
In the present we introduce a concept of doubly stochastic
quadratic operator. We prove necessary and sufficient conditions
for doubly stochasticity of operator. Besides, we prove that the
set of all doubly stochastic operators forms convex polytope.
Finally, we study analogue of Birkhoff's theorem for the class of
doubly stochastic operators.
{\it
Mathematics Subject Classification}: 15A51, 15A63, 46T99, 46A55.\\
{\it Key words}: Quadratic stochastic operator, doubly stochastic
operator, extremal point.

\end{abstract}

\maketitle

\section{Introduction}

Throughout the paper we will consider the simplex
$$S^{m-1}=\{x=(x_{1},
x_{2},\cdots
 , x_{m})\in R^{m}: x_{i}\geq 0 ,\ \forall i=\overline{1,m} \ \ \ \sum\limits_{i=1}^{m}x_{i}=1 \}.$$

For any $x=(x_{1},x_{2}\cdots x_{m})\in S^{m-1}$ due to \cite{Ma},
we define $x_{\downarrow}=(x_{[1]},x_{[2]},\cdots x_{[m]}),$ here
$x_{[1]}\geq x_{[2]}\geq \cdots \geq x_{[m]}$- nonincreasing
rearrangement of $x$ The point $x_{\downarrow}$ is called {\it
rearrangement} of $x$ by nonincreasing. Recall that for two
elements $x,y$ taken from the simplex $S^{m-1}$ we say that an
element $x$ {\it majorized} by $y$ (or $y$ {\it majorates} $x$),
and write $x\prec y$(or $y\succ x$) if the following hold
$$ \sum\limits_{i=1}^{k}x_{[i]}\leq \sum\limits_{i=1}^{k}y_{[i]}, for \ any \ k=\overline{1,m-1}.$$

Note  that such a term and notation was introduced by Hardy,
Littlewood and Polya in \cite{HLP}. It is easy to see that for any
$x\in S^{m-1}$ we have
$$(\frac{1}{m},\frac{1}{m},\cdots ,\frac{1}{m})\prec x\prec (1,0,\cdots 0).$$

A matrix $P=(p_{ij})_{i,j=\overline{1,m}}$ is called {\it doubly
stochastic(sometimes bistochastic)}, if
$$p_{ij}\geq 0, \ \ \forall i,j=\overline {1,m}$$
$$\sum\limits_{i=1}^mp_{ij}=1, \ \forall j=\overline {1,m}  \ \
 \ \sum\limits_{j=1}^mp_{ij}=1, \ \forall i=\overline {1,m}.$$

It is known \cite{Ma} that doubly stochasticity of a matrix $P$ is
equivalent to $Px\prec x$ for all $x\in S^{m-1}.$ It is clear that
the set of all doubly stochastic matrices is convex and compact.
Therefore, it was a problem due to Birkhoff concerning description
of the set of all extremal points of such a set. A solution of
that problem was given in \cite{bi}, and it states that the
extremal points consist of only permutations matrices.

The paper is devoted to the same problem mentioned above, but in a
class of nonlinear operators. To state the problem let us recall
some notions.

An operator $V:S^{m-1}\rightarrow S^{m-1}$ given by
$$(Vx)_{k}=\sum\limits_{i,j=1}^{m}p_{ij,k}x_{i}x_{j}.\eqno (1)$$
is called {\it quadratic stochastic operator (q.s.o.)}, here the
coefficients $p_{ij,k}$ satisfy the following conditions
$$p_{ij,k}=p_{ji,k}\geq 0 , \ \ \sum\limits_{k=1}^mp_{ij,k}=1.$$
One can see that q.s.o. is well defined, i.e. it maps the simplex
into itself. Note that such operators arise in many models of
physics , biology and so on. A lot of papers were devoted to
investigations of such operators (see, for example
\cite{Ber,Ga1,Lyu,M1,M2,M3}). We mention that in those papers
authors studied a central problem in the theory of q.s.o., namely
limit behavior of the trajectory of q.s.o.

We say that a q.s.o. (1) is called {\it doubly stochastic} (here
we saved the same terminology as above) if
$$Vx\prec x\eqno (2)$$ for
all $x\in S^{m-1}.$\footnote{Note that in mathematical economics
such an operator is usually called the operator of prosperity.}

The main object of the paper is the set of all doubly stochastic
quadratic operators (d.s.q.o.). The goal is to study necessary and
sufficient conditions for doubly stochasticity of quadratic
stochastic operators and extremal points of such the set of all
doubly stochastic quadratic operators. To do it, first in sec. 2
we study some properties of doubly stochastic q.s.o. Then we give
necessary conditions for doubly stochasticity. In section 3 we,
somehow, describe the set of doubly doubly stochastic operators.
Finally, in section 4 first of all we study a sufficient condition
for q.s.o. to be doubly stochastic and then study extreme points
of the set of such operators. Note that a part of the results was
announced in \cite{Ga2},\cite{Sh}.

\section{ Necessary conditions for doubly stochasticity of operators}

This section is devoted to some properties of the doubly
stochastic q.s.o.

\begin{thm}\label{1} Let $V:S^{m-1}\rightarrow S^{m-1}-$
be a d.s.q.o.. Then the coefficients $p_{ij,k}$ satisfy the
following conditions

$$a)   \sum\limits_{i,j=1}^mp_{ij,k}=m, \ \ \ \forall k=\overline{1,m}.\eqno (3)$$
$$b)   \sum\limits_{j=1}^mp_{ij,k}\ge\frac12, \ \ \ \forall i,k=\overline{1,m}.\eqno (4)$$
$$c)   \sum\limits_{i,j\in \alpha}p_{ij,k}\le |\alpha|, \ \ \forall \alpha\subset I, \ \ k=\overline{1,m}.\eqno (5)$$
here $I=\{1,2,\cdots,m\},$  $|\alpha|-$ cardinality of  $\alpha.$
\end{thm}
\begin{pf}

a)  Let $C=(\frac{1}{m},\frac{1}{m},\cdots
 ,\frac{1}{m}).$ From the definition we have $V(C)\prec C$, on the other hand $C\prec x, \ \forall x\in S^{m-1}.$
 Therefore
 $V(C)_{\downarrow}=
 C_{\downarrow}$ and since $C_{\downarrow}= C$ one gets $V(C)=C$ which implies (3).

c) Let $Vx\prec x.$ Then according to [Ma] there exists a doubly
stochastic matrix $P(x)=(p_{ij}(x))$ (which depends on $x$) such
that $Vx=P(x)x$. Put $x^0(\alpha)=(x_1^0,x_2^0,\cdots,x_m^0),$
here

$$\left\{
\begin{array}{ll}
x^0_i=\frac1{|\alpha|}, i\in\alpha\\
x^0_i=0, i\notin\alpha,\\
\end{array}
\right.
$$
where $\alpha$ is an arbitrary subset of $I$.

It is evident that $|\alpha|\leq m,$ therefore $x^{0}\in S^{m-1}$
and we have

$$(Vx^0)_k=\sum\limits_{i,j\in \alpha}^{m}p_{ij,k}\frac1{|\alpha|^2}=
\sum\limits_{i=1}^{|\alpha|}p_{ij,k}(x^0)\frac1{|\alpha|}\leq
\sum\limits_{i=1}^mp_{ij,k}(x^0)\frac1{|\alpha|}=\frac1{|\alpha|}.$$
Since $\alpha$ is an arbitrary set ,then we infer that
$$ \sum\limits_{i,j\in \alpha}p_{ij,k}\leq |\alpha|, \ \forall \alpha\subset I, \ \forall
 k=\overline{1,m}. $$

Using c) we are going to prove b).

b)  We will prove $ \sum\limits_{j=1}^mp_{i_0j,k}\ge\frac12$  for
some fixed $i_0.$

 From the equality (see (3))

$$m=\sum\limits_{i,j=1}^mp_{ij,k}=\sum\limits_{i,j=1,\ i,j\neq i_0 }^mp_{ij,k}+p_{i_0i_0,k}+2\sum\limits_{j=1,\ j\neq i_0}^mp_{i_0j,k}$$

and
$$\sum\limits_{i,j=1,\ i,j\neq i_0}^mp_{ij,k}\le m-1$$

we obtain
$$p_{i_0i_0,k}+2\sum\limits_{j=1,\ j\neq i_0}^mp_{i_0j,k}\ge 1$$

or

$$2\sum\limits_{j=1}^mp_{i_0j,k}\ge 1+p_{i_0i_0,k}\ge 1.$$
 Therefore, b) is satisfied.
\end{pf}

\section{Description of the class of doubly stochastic operators}

In this section  we are going to describe a class of doubly
stochastic operators.

Recall that a matrix $T=(t_{ij}), \ i,j=\overline{1,m}$ is said to
be {\it stochastic} if  $t_{ij}\ge 0$ and
$\sum\limits_{j=1}^mt_{ij}= 1.$

Let $A=(a_{ij}), \ i,j=\overline{1,m}$ be a symmetric matrix, with
$a_{ij}\ge 0.$

Consider the following equation with respect to $T$:

$$A=\frac12(T+T'),\eqno (6)$$
here $T'$ is the transposed matrix.

Below we are going to study conditions for solvability of equation
(6) in the class of all stochastic matrices.

Let  $G_m$ be the group of permutations of $m$ elements. For $g\in
G_m$ by  $A_g$ we denote a matrix $A_g=(a_{g(i)g(j)}), \ \
i,j=\overline{1,m},$ which is called {\it row and column
permutation} of
 $A.$

The following assertions are evident:

$i)$$$(A_g)_{g^{-1}}=A$$ for any $g\in G_m$

$ii)$ If $A$ symmetric, then  $A_g$ is also symmetric for any
$g\in G_m$

$iii)$ If $A$ stochastic, then  $A_g$ is also stochastic for any
$g\in G_m$

$iv)$ $$(A+B)_g=A_g+B_g ; \ \ (\lambda A)_g=\lambda A_g.$$

From  $i)-iv)$ we conclude that if  $T$ is a solution of (6), then
$T_g$ is a solution of the equation $$A_g=\frac12(T+T')$$

At first we will study the following set

$$\textbf{U}_1=\{A=(a_{ij}):a_{ij}=a_{ji}\ge 0, \ \ \sum\limits_{i,j\in\alpha}a_{ij}\le |\alpha|,  \ \ \sum\limits_{i,j\in I}a_{ij}=m\}.$$

It is easy to see that the above set is convex and compact. Below
we will study its extremal points.Let us shortly recall some
necessary notations.Let $A\subset X$ be convex set and $X$ be some
vector space. A point $x\in A$ is called extremal, if from
$2x=x_1+x_2,$ where $x_1, \ x_2\in A$ and $x_1\neq x_2,$ it
follows that $x=x_1=x_2.$ The set of all extremal points of a
given set $A$ is denoted by $extrA$ The Krein-Milman theorem
asserts that any convex compact set on some topological vector
space is a completion of the convex hull of its extremal points
(see, for review \cite{R}).

\begin{thm}\label{22} If $A=(a_{ij})\in extr\textbf{U}_1$ then $a_{ii}=0\vee 1, \ \ a_{ij}=0\vee \frac12\vee 1$
\end{thm}

Here and henceforth $c=a\vee b$ means that $c$ is either $a$ or
$b.$

\begin{pf} We prove this by induction with respect to the
order of the set $\textbf{U}_1.$ Let $m=2.$

First we prove that $a_{11}=0\vee 1$ and $a_{22}=0\vee 1.$ Let
$0<a_{11}<1,$ then from $a_{11}+2a_{12}+a_{22}=2$ we get
$0<a_{12}<1.$ Let us consider the following matrices

$A_1=\left(%
\begin{array}{cc}
  a_{11}+2\varepsilon & a_{12}-\varepsilon \\
  a_{12}-\varepsilon & a_{22} \\
\end{array}%
\right)$ and

$A_2=\left(%
\begin{array}{cc}
  a_{11}-2\varepsilon & a_{12}+\varepsilon \\
  a_{12}+\varepsilon & a_{22} \\
\end{array}%
\right)$

Since $0<a_{11}<1$ and $0<a_{12}<1,$ then one can choose
$\epsilon$ such that $A_1,A_2\in \textbf{U}_1.$ Now we obtain that
$2A=A_1+A_2.$ Therefore $A\notin extr\textbf{U}_1.$ The last
contradiction shows that $a_{11}=0\vee1.$ By this analogy one can
prove that $a_{22}=0\vee 1$ Since $a_{11}=0\vee 1$ and
$a_{22}=0\vee 1$ it follows from $A\in \textbf{U}_1$ that
$a_{ij}=0\vee \frac12\vee 1.$ One can easily see that the number
of such matrices is 4. That is the following matrices:

$\left(%
\begin{array}{cc}
  1 & 0 \\
  0 & 1 \\
\end{array}%
\right)$ $\left(%
\begin{array}{cc}
  0 & 1 \\
  1 & 0 \\
\end{array}%
\right)$ $\left(%
\begin{array}{cc}
  1 & \frac12 \\
  \frac12 & 0 \\
\end{array}%
\right)$ $\left(%
\begin{array}{cc}
  0 & \frac12 \\
  \frac12 & 1 \\
\end{array}%
\right)$

And all of them are extremal. Thus, for $A\in extr\textbf{U}_1$ it
is necessary and sufficient to be $a_{ii}=0\vee 1, \ \
a_{ij}=0\vee \frac12\vee 1.$ However for $m\ge 3$ it is not true
at all.

Let us suppose that the conditions of the theorem is valid for all
matrices of order less than $m$ and prove it for matrix of order
$m.$

First we prove that if $A\in extr\textbf{U}_1$ then $a_{ii}=0\vee
1.$ Let us assume that it is not so, that is there exists some
diagonal entry, which neither $0$ nor $1.$ Without loss of
generality we may assume that $0<a_{11}<1$.

We call a set $\alpha\in I$ is said to be \textit{saturated}, if
$\sum\limits_{i,j\in\alpha}a_{ij}=|\alpha|.$ Corresponding minor
to the saturated set we call \textit{saturated minor}. For
instance, $I$ itself is saturated, but it is a proper set.
Further, we mean only non proper saturated sets.

If there is a saturated set, such that $1\in \alpha$ then by the
assumption of the induction we infer that $a_{11}=0\vee 1.$ Now
remains the case when $1\notin\alpha$ for any saturated set
$\alpha.$

Now we will prove that if $0<a_{11}<1$ then $a_{jj}=0\vee 1$ for
all $j.$ Indeed, if there is $j_0$ such that $0<a_{j_0j_0}<1,$
then we put

$$A'=\{a_{11}'=a_{11}+\varepsilon, \ a_{j_0j_0}'=a_{j_0j_0}-\varepsilon, \ a_{ij}'=a_{ij}\ for \ all\ other\ values\ of\ i,j \}$$
and
$$A''=\{a_{11}''=a_{11}-\varepsilon, \ a_{j_0j_0}''=a_{j_0j_0}+\varepsilon, \ a_{ij}''=a_{ij}\ for\ all\ other\ values\ of\ i,j \}$$

If $j_0$ is contained in some saturated set, then from the
assumption of the induction it directly follows that
$a_{j_0j_0}=0\vee 1.$ If it is not so, then $A',A''\in
\textbf{U}_1$ and consequently $2A=A'+A'',$ that is $A\notin
\textbf{U}_1$, which is a contradiction.

So if $0<a_{11}<1,$ then $a_{jj}=0\vee 1$ for all $j.$
Furthermore, since $I$ is a saturated set, then there is
$a_{i_0j_0}\neq 0\vee\frac12\vee1(i_0\neq j_0).$ Now we put

$$A'=\{a_{11}'=a_{11}+2\varepsilon, \ a_{i_0j_0}'=a_{j_0j_0}-\varepsilon, \ a_{ij}'=a_{ij}\ for\ all\ other\ values\ of\ i,j \}$$
and
$$A''=\{a_{11}''=a_{11}-2\varepsilon, \ a_{i_0j_0}''=a_{j_0j_0}+\varepsilon, \ a_{ij}''=a_{ij} for\ all\ other\ values\ of\ i,j \}$$

If $\{i_0,j_0\}$ is not contained in any saturated set, then
$A',A''\in \textbf{U}_1$ and $2A=A'+A''.$

If $\{i_0,j_0\}$ is contained in some saturated set. Then by the
assumption of the induction we get $a_{i_0j_0}0=
0\vee\frac12\vee1(i_0\neq j_0),$ which contradicts to
$a_{i_0j_0}\neq 0\vee\frac12\vee1(i_0\neq j_0).$

Thus, all cases come to a contradiction. Therefore, $a_{11}=0\vee
1.$ Now we have to prove that if $A\in extr\textbf{U}_1$ then
$a_{ij}=0\vee\frac12\vee 1(i\neq j).$

Let  $a_{i_0j_0}\neq 0\vee\frac12\vee1(i_0\neq j_0).$ If there is
a saturated set $\alpha$ such that $\{i_0,j_0\}\subset \alpha$
then by the assumption it follows that $a_{ij}=0\vee\frac12\vee
1(i\neq j).$ Since $I$ is saturated, one can choose $\{i_1,j_1\}$
such that $a_{i_1j_1}\neq 0\vee\frac12\vee1(i_0\neq j_0).$ Also it
follows that $\{i_1,j_1\}$  is not contained in any saturated set.

Let us put

$$A'=\{a_{i_0j_0}'=a_{i_0j_0}+\varepsilon, \ a_{i_1j_1}'=a_{i_1j_1}-\varepsilon, \ a_{ij}'=a_{ij}\ for\ all\ other\ values\ of\ i,j \}$$
and
$$A''=\{a_{i_0j_0}''=a_{i_0j_0}-\varepsilon, \ a_{i_1j_1}''=a_{i_1j_1}+\varepsilon, \ a_{ij}''=a_{ij} for\ all\ other\ values\ of\ i,j \}$$

The sets $\{i_0,j_0\}$ and $\{i_1,j_1\}$ is not contained in any
saturated set. Therefore $A',A''\in \textbf{U}_1$ and $2A=A'+A''.$
\end{pf}

\begin{cor}\label{20}
If $m=3,$ then $A\in extr\textbf{U}_1$ if and only if
$a_{ii}=0\vee 1, \ \ a_{ij}=0\vee \frac12\vee 1$ and $A\neq M.$
Here $M=\left(%
\begin{array}{ccc}
  0 & \frac12 & \frac12  \\
  \frac12  & 0 & \frac12  \\
  \frac12  & \frac12 & 0 \\
\end{array}%
\right)$

Moreover, $|extr\textbf{U}_1|=25$
\end{cor}

\begin{cor}\label{21} If $A\in extr\textbf{U}_1$ then either $A$ is
stochastic or has a saturated minor of order $m-1.$
\end{cor}

\begin{pf} Let $A\in extr\textbf{U}_1.$  Let us assume that the matrix $A$ has no saturated minors of
order $m-1.$ Then we have $\sum\limits_{i,j\in
\alpha}a_{ij}<|\alpha|$ for any $\alpha$ such that $|\alpha|=m-1.$
By taking into account $\sum\limits_{i,j=1}^ma_{ij}=m$, we obtain

$$a_{ii}+\sum\limits_{j\neq i}a_{ij}> 1.$$

By theorem \label{22}it follows that $a_{ii}=0\vee 1, \ \
a_{ij}=0\vee \frac12\vee 1.$ Therefore an expression
$a_{ii}+\sum\limits_{j\neq i}a_{ij}$ is an integer. Therefore

$$a_{ii}+\sum\limits_{j\neq i}a_{ij}\ge 2.$$

By taking summation in the above inequality from $i=1$ to $i=m$ we
find
$$\sum\limits_{i=1}^m(a_{ii}+\sum\limits_{j\neq i}a_{ij})=2\sum\limits_{i,j=1}^ma_{ij}-\sum\limits_{i=1}^ma_{ii}\ge 2m.$$

Hence
$$\sum\limits_{i=1}^ma_{ii}\le 0.$$
and $a_{ii}=0$ for all $i=\overline{1,m}.$ Further,
$a_{ii}+\sum\limits_{j\neq i}a_{ij}\ge 2$ and $a_{ii}=0$ imply
$\sum\limits_{j=1}^ma_{ij}\ge 1$ for all $i=\overline{1,m}.$

Since $\sum\limits_{i,j=1}^ma_{ij}=m$ then
$\sum\limits_{j=1}^ma_{ij}=1.$ The last means that $A$ is
stochastic.
\end{pf}

\begin{thm}\label{2} Let $A=(a_{ij})-$ be a symmetric and nonnegative matrix.
For the existence of a stochastic matrix $T=(t_{ij})$ satisfying
the equation (6) necessary and sufficient condition is
$$\sum\limits_{i,j\in\alpha}a_{ij}\le |\alpha|$$
for all $\alpha\subset I=\{1,2,\cdots,m\}.$
\end{thm}

\begin{pf}

\texttt{Necessity}. If $T=(t_{ij})$ is a stochastic matrix and (6)
is satisfied, then

$$\sum\limits_{i,j\in\alpha}a_{ij}=\frac12(\sum\limits_{i,j\in\alpha}t_{ij}+\sum\limits_{i,j\in\alpha}t_{ji})=
\sum\limits_{i,j\in\alpha}t_{ij}= \sum\limits_{i\in\alpha}
\sum\limits_{j\in\alpha}t_{ij}\le
\sum\limits_{i,j\in\alpha}1=|\alpha|$$

\texttt{Sufficiency.} Let $A\in \textbf{U}_1.$ We have to show
that existence of a stochastic matrix for which (6) is satisfied.

First, we prove it for extremal points of $\textbf{U}_1.$ Now let
$A\in extr\textbf{U}_1.$ Then from corollary \ref{21} it follows
that either $A$ is stochastic or has a saturated minor of order
$m-1.$ If $A$ is stochastic, then we take $T=A$ since $A$ is
symmetric then (6) is satisfied.

Now, let us suppose that $A$ is not stochastic, then it follows
that it has a saturated minor order $m-1.$ This case we prove by
induction relatively order of the matrix. Let $m=2$.
 In this case we have the following matrices
 $A_1=\left(%
\begin{array}{cc}
  1 & 0 \\
  0 & 1 \\
\end{array}%
\right)$,
$A_2=\left(%
\begin{array}{cc}
  0 & 1 \\
  1 & 0 \\
\end{array}%
\right)$,
$A_3=\left(%
\begin{array}{cc}
  1 & \frac12 \\
  \frac12 & 0 \\
\end{array}%
\right)$,
$A_4=\left(%
\begin{array}{cc}
  0 & \frac12 \\
  \frac12 & 1 \\
\end{array}%
\right)$.

For $A_1$ or $A_2$ we take $T=A_1$ and $T=A_2$, respectively.

For $A_3$ or $A_4$ we take

$T=\left(%
\begin{array}{cc}
  1 & 0 \\
  1 & 0 \\
\end{array}%
\right)$ or
$T=\left(%
\begin{array}{cc}
  0 & 1 \\
  0 & 1 \\
\end{array}%
\right)$, respectively.

Let us assume that assumption of the theorem is true for all
matrices of order less that $m.$

Since $A$ has a saturated set of order $m-1,$ without loss of
generality we may assume that this minor is
$(a_{ij})_{i,j=\overline{1,m-1}}.$ From the assumption of the
induction there is a stochastic matrix for this minor which
satisfies (6).Let it be $V=(v_{ij})_{i,j}=\overline{1,m-1}.$

Since $\sum\limits_{i,j=1}^{m-1}a_{ij}=m-1$ and
$\sum\limits_{i,j=1}^ma_{ij}=m$ it follows that
$a_{mm}+2\sum\limits_{i=1}^{m-1}a_{im}=1.$ From theorem \label{22}
it follows either $a_{mm}=1$ and $a_{im}=0$ for
$i=\overline{1,m-1}$ or $a_{ii}=0,$  $a_{i_0m}=\frac12$ $a_{im}=0$
for $i=\overline{1,m-1}\setminus\{i_0\}$

For the first case we define $T$ by

$$T\{(t_{ij}): t_{ij}=v_{ij},\  i,j=\overline{1,m-1} \  t_{mm}=1,\ t_{im}=t_{mi}=0,\ \i=\overline{1,m-1}\}$$
and it is easy to see that above matrix is stochastic and
satisfies (6).

For the second case we define $T$  by

$$T\{(t_{ij}): t_{ij}=v_{ij},\  i,j=\overline{1,m-1} \  t_{mi_0}=1,\ t_{im}=t_{mi}=0,\ \i=\overline{1,m}\}$$

and it is easy to see that above matrix is stochastic and
satisfies (6). So, for extremal points of $\textbf{U}_1$ the
theorem has been proved.

Now let $A\in \textbf{U}_1$ From Krein-Milman(see \cite{R}) it
follows that $A$ can be produced as the convex hull of its
extremal points.

Let $$A=\sum\limits_{i=1}^k\lambda_{k}A_{k}.$$

Here $0\le \lambda\le 1,\ \sum\limits_{i=1}^k\lambda_{k}=1$ and
$A_{i}\in extr\textbf{U}_1, \ i=\overline{1,k}.$

Let for matrices $A_{i}$ correspond stochastic matrices $T_{i}$
that satisfy (6). Then
$$T=\sum\limits_{i=1}^k\lambda_{k}T_{k}$$ is stochastic and satisfies (6).

\end{pf}

\textbf{Remark.} It should be mentioned that stochastic matrix $T$
(which is a solution of (6)) exists but not unique.

\textbf{Example.} For the matrix
$$A=\left(%
\begin{array}{ccc}
  0,1 & 0,3 & 0,4 \\
  0,3 & 0,1 & 0,5 \\
  0,4 & 0,5 & 0,4 \\
\end{array}%
\right)$$

the set of the solutions of $A=\frac12(T+T')$ in the class of
stochastic matrices is

$$T_{\alpha}=\left(%
\begin{array}{ccc}
  0,1 & \alpha & 0,9-\alpha \\
  0,6-\alpha & 0,1 & 0,3+\alpha \\
  \alpha-0,1 & 0,7-\alpha & 0,4 \\
\end{array}%
\right)$$

here $\alpha\in [0,1;0,6].$\\

Recall that a set $A$ in $R^k$ is said to be a {\it polytope}, if
it is nonempty, bounded and consists of intersection of a finite
number of semispaces.
\begin{prop}
The set of all solutions of the equation
$$A=\frac12(T+T')$$
 with respect to $T$ forms a polytope.
\end{prop}
\begin{pf}
Let us rewrite the solutions of the equation in the following
form:
$$0\le t_{ij}\le 1,$$
$$t_{ij}+t_{ji}\le a_{ij}, \ \ t_{ij}+t_{ji}\ge a_{ij}$$
$$\sum\limits_{j=1}^mt_{ij}\le 1,\ \ \sum\limits_{j=1}^mt_{ij}\ge 1.$$

From above given relations one can see that the set of all
solutions of the equation is nonempty, bounded and it is
intersection of semispaces.

\end{pf}

From (6) it follows that $(Ax,x)=(Tx,x),$ here $(\cdot,\cdot)$ is
the inner product in $R^m.$ Indeed

$$(Ax,x)=\frac12\left((T+T')x,x\right)=\frac12\left[(Tx,x)+(T'x,x)\right]=(Tx,x).\eqno (7)$$

Now let's return to the study of d.s.q.o. According to the theorem
2.1., if $V:S^{m-1}\rightarrow S^{m-1}-$ d.s.q.o., then

$$   \sum\limits_{i,j\in \alpha}p_{ij,k}\le |\alpha|, \ \ \forall
\alpha\in I, \ \ k=\overline{1,m}$$  moreover, for $\alpha=I$ an
equality is held.

By putting
$$A_k=(p_{ij,k})_{i,j=\overline{1,m}},$$
and using the theorem\label{2} one can find a stochastic matrix
$T_k$ such that

$$A_k=\frac12(T_k+T'_k), \ \ k=\overline{1,m}.$$

From (7) one gets
$$\sum\limits_{i,j=1}^mp_{ij,k}x_ix_j=(A_kx,x)=(T_kx,x).$$
Hence, if $V$ is a  d.s.q.o., then there are stochastic matrices
$T_1,T_2,\cdots,T_m$ such that
$$Vx=\left((T_1x,x),(T_2x,x),\cdots,(T_mx,x)\right).$$

Now we prove the following inequality:

$$\min_{1\le i\le m}x_i\le (Tx,x)\le \max_{1\le i\le m}x_i, \ \ \forall x\in S^{m-1},$$
or in accepted notations:

$$x_{[m]}\le (Tx,x)\le x_{[1]}, \ \ \forall x\in S^{m-1},\eqno (8)$$
here $T-$ is a stochastic matrix.

Indeed, if $t_{ij}\ge 0$ and $\sum\limits_{j=1}^mt_{ij}=1,$ then
$$x_{[m]}\le \sum\limits_{j=1}^mt_{ij}x_{j}\le x_{[1]}$$
for all $x\in R^m.$ In particular,for $x\in S^{m-1}$ we get

$$(Tx,x)=\sum\limits_{i,j=1}^mt_{ij}x_ix_j=\sum\limits_{i=1}^mx_i\sum\limits_{j=1}t_{ij}x_j.$$

By $x_i\ge 0, \ \ \sum\limits_{i=1}^mx_i=1$ and above inequality
we get the required inequality (8).

\begin{thm}\label{4} Let $A=(a_{ij})$  be a nonnegative symmetric matrix.
Then for the fulfillment of the inequality
$$x_{[m]}\le (Ax,x)\le x_{[1]} \ \ \forall x\in S^{m-1},$$
it is necessary and sufficient that

$$   \sum\limits_{i,j\in \alpha}a_{ij}\le |\alpha| \ \ \forall
\alpha\in I, \ \ k=\overline{1,m}$$  moreover, for $\alpha=I$ an
equality is held.
\end{thm}
\begin{pf} Sufficiency straightforwardly follows from Theorem 3.4 and the inequality (8).

 {\it Necessity}. Consider a set $\alpha\subset I$ and put $x^0=(x_1^0,x_2^0,\cdots,x_m^0),$  with

 $$x_i^0=\frac1{|\alpha|}, \ \ i\in \alpha, \ \ \ x_i^0=0, \ \ \ i\notin \alpha.$$
Then  for $|\alpha|<m$ we get $x_{[m]}=0, \
x_{[1]}=\frac1{|\alpha|}$. Therefore, one has

$$0\le (Ax,x)=\sum\limits_{i,j=1}^ma_{ij}x_ix_j=\frac1{|\alpha|^2}\sum\limits_{i,j\in \alpha}a_{ij}\le \frac1{|\alpha|}$$
hence
$$\sum\limits_{i,j\in \alpha}a_{ij}\le |\alpha|$$

If  $|\alpha|=m,$ then $x_{[m]}=x_{[1]}=\frac1m$. Consequently
$$\frac1m\le \frac1{m^2}\sum\limits_{i,j\in \alpha}a_{ij}\le \frac1m,$$
 that is  $\sum\limits_{i,j=1}^ma_{ij}=m.$
\end{pf}
 By Theorems 3.4 and 3.5 one can deduce that an inequality

$$x_{[m]}\le (Ax,x)\le x_{[1]} \ \ \forall x\in S^{m-1}\eqno (9)$$
here $A=(a_{ij})$ is a nonnegative matrix, is equivalent to the
existence of stochastic matrix $T=(t_{ij})$,
 such that
 $$A+A'=T+T'.$$

Since we are going to study sufficient conditions for doubly
stochasticity of q.s.o. from the proved theorem one arises  the
following

\begin{prob} For which symmetric matrices $A$ the following inequality holds

$$x_{[m]}+x_{[m-1]}+\cdots+x_{[m-k+1]}\le (Ax,x)\le x_{[1]}+\cdots+x_{[k]},\eqno (10)$$
for all $x\in S^{m-1},$ here $k-$ natural number $1\le k\le m ?$
\end{prob}

For $k=1$  an answer is given by theorem 3.3. To solve the problem
we give some notations:

$$\textbf{T}_k=\{T=(t_{ij}), \ i,j=\overline{1,m}: 0\le t_{ij}\le 1, \sum\limits_{j=1}^mt_{ij}=k\}, \ \ 1\le k\le m.$$

$$\textbf{U}_k=\{A=(a_{ij}):a_{ij}=a_{ji}, \ \ A=\frac12(T+T'), \ T\in \textbf{T}_{k}  \}, \ \ 1\le k\le m.$$

The set $\textbf{U}_k$ is called a \textit{symmetrization} of
$\textbf{T}_k.$ Evidently, $\textbf{T}_1-$ the set of all
stochastic matrices.

\begin{thm}\label{5} For $A\in \textbf{U}_k$ the inequality (10) holds.
\end{thm}
\begin{pf}Let
$x_{\downarrow}=(x_{[1]},\cdots,x_{[m]})-$ be a  nonincreasing
rearrangement of $x$ and $\lambda_1,\cdots,\lambda_m$ be an
arbitrary numbers such that $0\le \lambda_i\le 1,$
$\sum\limits_{i=1}^m\lambda_i=k.$ Consider the following sum
$$\lambda_1x_{[1]}+\lambda_2x_{[2]}+\cdots+\lambda_mx_{[m]}.\eqno (11)$$
Let, $i<j.$ Replacement of the coefficients
$\lambda_1,\cdots,\lambda_i,\cdots,\lambda_j,\cdots,\lambda_m$
into
$\lambda_1,\cdots,\lambda_i+\varepsilon,\cdots,\lambda_j-\varepsilon,\cdots,\lambda_m$
we will call \textit{backward shift}, if $\varepsilon>0$ and
\textit{forward shift} if $\varepsilon<0.$ Shift is called
\textit{admissible}, if it preserves the condition $0\le
\lambda'_i,$ here $\lambda'_i$ are the coefficients, obtained in
the result of shift.

Since $x_{[i]}\ge x_{[j]},$ then it is clear that under the
admissible backward shifts, the sum (11) does not decrease, and
respectively, does not increase under the admissible forward
shifts.

It is easy to see that an admissible backward shift is possible up
to obtaining the collection $1,1,\cdots,1,0,0,\cdots,0,$ and to
the right up to obtaining the collection
$0,0,\cdots,0,1,1,\cdots,1,$ here in both cases the number of ones
is $k.$ That's why

$$x_{[m-k+1]}+x_{[m-k+2]}+\cdots+x_{[m]}\le \lambda_1x_{[1]}+\lambda_2x_{[2]}+\cdots+\lambda_mx_{[m]}\le x_{[1]}+\cdots+x_{[k]} \eqno (12)$$
 for an arbitrary $x\in R^m$ and $0\le \lambda_i\le 1, \ \ \sum\limits_{i=1}^m\lambda_i=k.$

Now let $x\in S^{m-1}$ and $A\in \textbf{U}_k.$

Choose $T\in \textbf{T}_k$,such that

$$A=\frac12(T+T')$$

Then
$$(Tx,x)=\sum\limits_{i,j=1}^mt_{ij}x_ix_j=\sum\limits_{i=1}^mx_i(\sum\limits_{j=1}^mt_{ij}x_j)$$
Since $0\le t_{ij}\le 1$ and $\sum\limits_{j=1}^mt_{ij}=k,$ from
(12) one gets the following

$$\sum\limits_{i=m-k+1}^mx_{[i]}\le \sum\limits_{j=1}^mt_{ij}\le \sum\limits_{i=1}^kx_{[i]} $$

According to $x_i\ge 0, \ \sum\limits_{i=1}^mx_i=1$ we obtain

$$\sum\limits_{i=m-k+1}^mx_{[i]}\le (Tx,x) \le \sum\limits_{i=1}^kx_{[i]} $$

or
$$\sum\limits_{i=m-k+1}^mx_{[i]}\le (Ax,x) \le \sum\limits_{i=1}^kx_{[i]} $$

So if $A\in \textbf{U}_k$ then (10) holds.
\end{pf}

Let $\textbf{B}$ be the set of all d.s.q.o. By putting
$A_k=\{p_{ij,k}\}_{(i,j=\overline{1,m})}$ we rewrite an operator
$V$ in the following form

$$V=(A_{1}|A_{2}|\cdots |A_{m}).\eqno (13)$$

Then theorems 2.1, 3.4, 3.5 and 3.6 imply that  for $V\in
\textbf{B}$ the conditions
$$A_k\in \textbf{U}_1, \ \ k=\overline{1,m} \eqno (14)$$
are necessary and the conditions

$$\forall \alpha\subset I, \ \sum\limits_{k\in \alpha}A_k\in \textbf{U}_{|\alpha|}\eqno (15)$$

are sufficient.

Whether the conditions (14) is  sufficient to be $V\in
\textbf{B}?$ We give a positive answer for the question in next
section.

Now we show some simple properties of the $\textbf{U}_k$ which
will be helpful. Let $E$ be the matrix of $m\times m$ with all
entries equal to unit.

\begin{thm}\label{6}
The following assertions hold

$$i)  \ A\in \textbf{U}_k\Leftrightarrow E-A\in \textbf{U}_{m-k};$$
$$ii) \ \textbf{U}_k\cap \textbf{U}_l=\emptyset, \ \ k\neq l;$$
$$iii) \ \textbf{U}_{m}=\{E\};$$
$$iv)\  A\in \textbf{U}_k\Rightarrow \frac{p}{k}A\in \textbf{U}_{p},
1\le p\le k;$$
$$v) \ \textbf{U}_k+\textbf{U}_l\supset \textbf{U}_{k+l}, \ \ k+l\le m .$$

\end{thm}
\begin{pf}
i) Let $1\le k\le m-1$ and $A\in U_k.$ From
$$A=\frac12(T+T')$$
here $T=(t_{ij}), \ 0\le t_{ij}\le 1, \
\sum\limits_{j=1}^mt_{ij}=k$ one gets
$$E-A=\frac12[(E-T)+(E-T)']$$
It is obvious that $E-T\in \textbf{T}_{m-k}$, therefore $E-A\in
\textbf{U}_{m-k}.$

ii) It follows from the definition of $\textbf{U}_k.$

iii) Since $A\in \textbf{U}_k$ then $0\le a_{ij}\le 1$ for all
$k=\overline{1,m}.$ If $A\in \textbf{U}_m$ then
$\sum\limits_{i,j=1}^ma_{ij}=m^2$ and therefore $a_{ij}=1$ that is
$\textbf{U}_{m}=\{E\}$

iv)Let  $A\in U_k$ and $A=\frac12(T+T').$ Then
$$\frac{p}{k}A=\frac12[\frac{p}{k}T+\frac{p}{k}T']$$
moreover $\frac{p}{k}T\in \textbf{T}_p.$ Consequently,
$\frac{p}{k}A\in \textbf{U}_p.$

v) Let $k+l\le m, A\in \textbf{U}_{k+l}$ and $A=\frac12(T+T')$ By
denoting

$$A_1=\frac{k}{2(k+l)}(T+T'), \ \ A_2=\frac{l}{2(k+l)}(T+T')$$ we
obtain $A=A_1+A_2,$ furthermore  $A_1\in \textbf{U}_k,
A_2\in\textbf{ U}_l.$ The last implies
$$\textbf{U}_k+\textbf{U}_l\supset \textbf{U}_{k+l}, \ \ k+l\le m
.$$ It should be mentioned that the last inclusion is strict for
all $k,l\ge 1$ Indeed, let $I$ be the identity matrix. Clearly
that $I\in \textbf{U}_1,$ but $I+I\notin\textbf{ U}_2.$
\end{pf}

\section{Extreme point of the set of doubly stochastic operators and Birkhoff's problem}

 According to the classic result of Birkhoff (see \cite{bi}) extreme
 points of the set of doubly stochastic matrices are permutation
 matrices, i.e. matrices having exactly one unit entry in each row and
 exactly one unit entry in each column, all other entries being
equal to zero. It is interesting to know an answer for the similar
problem about the set of all doubly stochastic nonlinear
operators. Therefore, to investigate such a problem it is better
first to start with the set of all d.s.q.o.'s, since such a set
contains as a subset of the set of all doubly stochastic matrices.
In this section we are going to describe extreme points of the set
of d.s.q.o.  At first, we prove necessary and sufficient
conditions for q.s.o.
 to be d.s.q.o.(it was proved in the last section but in this section
 we show that necessary and sufficient will coincide).

The set of all d.s.q.o. we denoted by $\textbf{B}$. Directly from
definition it follows that

$$V=(A_1|A_2|\cdots|A_m)\in \textbf{B}\Leftrightarrow V=(A_{\pi(1)}|A_{\pi(2)}|\cdots|A_{\pi(m))}\in \textbf{B}$$
here $\pi$ is an arbitrary permutation of the index set
$I=\{1,2,\cdots,m\}$

\begin{thm}\label{7}
Assume that $A_1\in\textbf{ U}_1.$ Then one can choose
$A_2,\cdots, A_m\in \textbf{U}_1$ such that
$V=(A_1|A_2|\cdots|A_m)\in \textbf{B}.$
\end{thm}
\begin{pf}
Take

$$A_2=A_3=\cdots=A_m=\frac1{m-1}(E-A_1)$$
here $E$ is the matrix with all entries equal to unit. According
to the theorem 3.8 we have
$$E-A_1\in \textbf{U}_{m-1}$$
and
$$\frac1{m-1}(E-A_1)\in \textbf{U}_1$$

Let a stochastic matrix $T_1\in \textbf{T}_1$ be a solution of the
equation
$$A_1=\frac12(T+T')$$
Now we put
$$T_2=T_3=\cdots=T_m=\frac1{m-1}(E-T_1)$$

we show that the sum of any $k$ of the matrices $T_1,T_2,\cdots
T_m$ belongs to $\textbf{T}_k.$ Indeed, let $T_1=\{t_{ij}\}$ such
that $0\le t_{ij}\le 1.$ Then
$$T_2=T_3=\cdots=T_m=\frac1{m-1}\{1-t_{ij}\}$$
From $0\le t_{ij}\le 1$ one  gets
$$0\le \frac{k}{m-1}(1-t_{ij})\le 1 \eqno (16)$$
and
$$0\le t_{ij}+\frac{k-1}{m-1}(1-t_{ij})\le 1. \eqno (17)$$

If the sum contains $T_1$ then our assertion follows from the
inequality (17), otherwise from the inequality (16). So the sum of
any $k$ matrices of $T_1,\cdots, T_m$ belongs to $\textbf{T}_k.$

Now we return to $A_1,\cdots, A_m$ and obtain that the sum of any
of them belongs to $\textbf{U}_k.$ Therefore we get that

$$V=(A_1|A_2|\cdots|A_m)\in \textbf{B}$$

\end{pf}

It is clear that if $T_1,\cdots,T_p\in \textbf{T}_1(p\le m)$ and
$\sum\limits_{i=1}^pT_i\in \textbf{T}_p$ then for any $k\le p$ the
sum of any $k$ matrices of $T_1,\cdots,T_p$ belongs to
$\textbf{T}_k.$ This implies that from $A_1,\cdots,A_p\in
\textbf{U}_1, \ \sum\limits_{i=1}^pA_i\in \textbf{U}_p$ it follows
that the sum of any $k(k\le p)$ matrices belongs to
$\textbf{U}_k.$ Using the last assertion we get the following
\begin{cor}\label{8}
If $A_1,\cdots,A_p\in \textbf{U}_1, \ \sum\limits_{i=1}^pA_i\in
\textbf{U}_p(p<m)$ then one can choose $A_{p+1},\cdots A_m\in
\textbf{U}_1$ such that
$$V=(A_1|A_2|\cdots|A_m)\in \textbf{B}$$
\end{cor}
So conditions (15) can be changed to conditions:
$$A_i\in\textbf{U}_1, \ \ i=\overline{1,m}, \ \ \sum\limits_{i=1}^mA_i=E \eqno (18)$$
since $\textbf{U}_m$ is $E.$ Therefore

\begin{cor}\label{9}
The conditions (18)are necessary and sufficient for q.s.o. to be
d.s.q.o.

\end{cor}

\begin{cor}\label{10}
The set of all d.s.q.o. forms a convex polytope.
\end{cor}
\begin{pf}
We recall that finite intersection  of nonempty, bounded and
closed semispaces is called convex polytope (see \cite{gr}). For
any q.s.o. there is a cubic matrix which can be embedded in the
space $R^{m^3}$.  One can embed it this cubic matrix in
$R^{\frac{m^2(m-1)}2}.$ Furthermore, each of the condition (18)
defines closed semispace in $R^{\frac{m^2(m-1)}2}$. If we denote
the set of all d.s.q.o. by $\textbf B$, from the corollary 4.3 it
follows that consists of intersection of these semispaces. It is
clear that this intersection is nonempty and bounded set in
$R^{\frac{m^2(m-1)}2},$ since it lies in positive ortant and in
hyperplane $\sum\limits_{i,j,k=1}^{m}p_{ij,k}=m^{2}.$ Due to
Grunbaum \cite{gr} we conclude that $\textbf{B}$ is convex
polytope.
\end{pf}

In \cite{Ga1} it was conjectured that an operator
$V=(A_{1}|A_{2}|\cdots|A_{m})$ is extremal point of $\textbf{B}$
if and only if $A_{i}\in extr U, \ \forall i=\overline{1,m}$. We
show that in general it is not true.

\begin{thm}\label{15} Let  $V=(A_{1}|\cdots |A_{m})\in extr \mathbf B$.
Then $V_{\pi}=(A_{\pi (1)}|\cdots |A_{\pi (m)})\in extrU$ for any
permutation $\pi$ of the index set $\{1,2,\cdots, m\}$.
\end{thm}
\begin{pf} Let $V_{\pi}\notin extr \textbf B$.Then $\exists
V',V''\in \mathbf B, \ V'\neq V''$, such that $2V_{\pi}=V'+V''$.
Let $V'=(A_{1}'|\cdots |A_{m}'), \ V''=(A_{1}''|\cdots |A_{m}'')$.
Then $(2A_{\pi (1)}-A_{1}'-A_{1}''|\cdots |2A_{\pi
(m)}-A_{m}'-A_{m}'')= \textbf{0}$. Since the matrices $2A_{\pi
(i)}-A_{i}'-A_{i}''$ are symmetric, then $2A_{\pi
(i)}=A_{i}'+A_{i}''$, therefore $2A_{i}=A_{\pi^{-1}
(i)}'+A_{\pi^{-1} (i)}''$. Since
$$\sum\limits_{i=1}^{m}A_{i}'=\sum\limits_{i=1}^{m}A_{i}''=E$$
then
$$\sum\limits_{i=1}^{m}A_{\pi^{-1} (i)}'= \sum\limits_{i=1}^{m}A_{\pi^{-1} (i)}''=E.$$
 Denote $W'=(A_{\pi^{-1} (1)}'|\cdots|A_{\pi^{-1} (m)}')$ and
$W''=(A_{\pi^{-1} (1)}''|\cdots|A_{\pi^{-1} (m)}'')$. Then we get
$W',W''\in \textbf B, \ W'\neq W''$ è $2V=W'+W''$ which
contradicts to  $V\notin extr\textbf B.$

\end{pf}

\begin{thm}\label{16} Let
$V=(A_{1}|\cdots|A_{m}).$ If for any permutation $\pi$ of any
$m-1$ elements of
 $\{1,2,\cdots m\}$ we have $A_{\pi(k)}\in extr\textbf{U}$ ,  then $V\in extr\textbf B.$
\end{thm}

\begin{pf}Using the above theorem it is enough to show that if
 $A_{1},A_{2},\cdots,A_{m-1}\in extr\textbf{U}$, then $V\in extr\textbf B$.
Let us assume that $V\notin extr \textbf B$. Then$\exists
V',V''\in \mathbf B ,\ \ V'\neq V''$ such that $2V=V'+V''$.Let
$V'=(A_{1}'|\cdots|A_{m}'), \ \ V''=(A_{1}''|\cdots|A_{m}'')$.
Then
$$(2A_{1}-A_{1}'-A_{1}''|\cdots |2A_{m}-A_{m}'-A_{m}'')=
\textbf{0}.$$ So $2A_{i}-A_{i}'-A_{i}''=\textsl{0}\, \ \forall
i=\overline{1,m}$. Since $A_{i}\in extrU,  \ \forall i\geq 2$,
then $A_{i}'=A_{i}'', \ \ \forall i\geq 2$. From the equality
$\sum\limits_{i=1}^{m}A_{i}'=\sum\limits_{i=1}^{m}A_{i}''=E$ we
get $A_{1}'=A_{1}''$. That's why $V'=V''$, which contradicts to
$V\notin extr\textbf B$.
\end{pf}

\begin{cor}\label{17} If $A_{i}\in extr \textbf{U}, \
\forall i=\overline{1,m}$, then $V\in extr \mathbf{B}$. However,
the converse case is not true
\end{cor}

\textbf{Example :} Consider an operator
$$\begin{array}{c}
  x_{1}'=x_{1}x_{2}+x_{1}x_{3}+x_{2}x_{3} \\
  x_{2}'=x_{2}^{2}+ x_{3}^{2}+x_{1}x_{3} \\
  x_{3}'=x_{1}^{2}+x_{1}x_{2}+x_{2}x_{3} \\
\end{array}$$ Corresponding matrices  have the following view $$\left(%
\begin{array}{ccc}
  0 & \frac{1}{2} & \frac{1}{2} \\
  \frac{1}{2} & 0 & \frac{1}{2} \\
  \frac{1}{2} & \frac{1}{2} & 0 \\
\end{array}%
\right)
\left(%
\begin{array}{ccc}
  0 & 0 & \frac{1}{2} \\
  0 & 1 & 0 \\
  \frac{1}{2} & 0 & 1 \\
\end{array}%
\right) \left(%
\begin{array}{ccc}
  1 & \frac{1}{2} & 0 \\
  \frac{1}{2} & 0 & \frac{1}{2} \\
  0 & \frac{1}{2} & 0 \\
\end{array}%
\right)$$

 It is easy to see that two of them are extremal.Then from the above theorem it
 follows that $V\in extr \textbf B$.

Whether the conditions of the theorem \label{16} are necessary and
sufficient? In two dimensional simplex the problem is solved
positively.

\begin{thm}\label{18} Let $V=(A_{1}|A_{2}|A_{3})\in \textbf B.$ $V\in
extr\textbf B$ if the only if at least 2 of the matrices
$A_{1},A_{2},A_{3}$ are extremal of $\textbf{U}_1$.
\end{thm}

\begin{pf}

For the proof of the given theorem we refer to the theorem 5 from
\cite{Ga2}, which says that if $V\in extr\mathbf B$, then entries
of the matrices  $A_{1},A_{2},A_{3}$ either $0$ or $\frac12$ or
$1.$ Therefore, if $V\in extr\mathbf B$, then either
$A_{1},A_{2},A_{3}\in extr\mathbf U$ or $A_{1},A_{2},A_{3}=M.$

The case $A_{1}=M,\ A_{2}=M, A_{3}=M$ is impossible because of
$A_1+A_2+A_3=E.$

It can be easily shown that if two of the matrices
$A_{1},A_{2},A_{3}$ is $M$ then $V$ is not extremal. Therefore we
can conclude that at least two of the matrices $A_{1},A_{2},A_{3}$
are extremal.
\end{pf}

\begin{cor}\label{19} For $m=3$ we have $|extr\mathbf B|=222.$
\begin{pf}
Let $V\in extr\mathbf B, \ \ V=(A_{1}|A_{2}|A_{3}).$ From the
theorem 4.8 it follows that $A_{2},A_{3}\in extr\mathbf U$ ñ up to
permutation. From $A_{1}+A_{2}+A_{3}=E$ one gets that either
$A_{1}\in extr\mathbf U$ or $A_{1}=M.$

Let $A_{1}\in extr\mathbf U.$ From  corollary 3.2 we know all
extreme points of $\mathbf U$ Therefore we can choose those
triples of extreme points, sum of which is $E.$ The number of such
triples is 31.For the case of  $A_{1}=M,$ the number of such
triples is 6.

Consequently,  $|extr\mathbf B|=37\times 3!=222.$
\end{pf}

\end{cor}

\end{document}